\documentclass{amsart}

\usepackage{amsmath,amsthm,amssymb,amscd,enumitem}
\usepackage{amsfonts}
\usepackage{rotating}
\usepackage{euscript}
\usepackage{pst-node}
\usepackage{epsfig,verbatim,mathtools}
\mathtoolsset{showonlyrefs}
\def\be{\begin{equation}}
\def\ee{\end{equation}}

\def\Aut{{\rm Aut}}

\def\C{{\mathbb C}} 
\def\f{\EuScript}
 
\def\P{{\mathbb P}}

\def\phi{{\varphi}}
 
\def\tt{\widetilde}
\def\deg{{\rm deg\,}}

\def\Ker{{\rm Ker\,}}

\def\GCD{{\rm GCD }}

\def\mod{{\rm mod\ }}

\def\bp{\begin{proposition}}
\def\ep{\end{proposition}}

\def\bt{\begin{theorem}}
\def\et{\end{theorem}}
\def\br{\begin{remark}}
\def\er{\end{remark}}
\def\be{\begin{equation}}
\def\bee{\begin{equation*}}
\def\l{\label}

\def\ee{\end{equation}}
\def\eee{\end{equation*}}
\def\bl{\begin{lemma}}
\def\el{\end{lemma}}
\def\bc{\begin{corollary}}
\def\ec{\end{corollary}}
\def\pr{\noindent{\it Proof. }}

\def\bd{\begin{definition}}
\def\ed{\end{definition}}
\def\t{\widetilde}
\def\h{\widehat}

\mathtoolsset{showonlyrefs}

\newtheorem{theorem}{Theorem}[section]
\newtheorem{lemma}[theorem]{Lemma}
\newtheorem{definition}[theorem]{Definition}
\newtheorem{corollary}[theorem]{Corollary}
\newtheorem{proposition}[theorem]{Proposition}

\theoremstyle{definition}

\theoremstyle{definition}
\newtheorem{remark}[theorem]{Remark}

\begin{document}

\title{Hurwitz existence problem and fiber products}
\author{Fedor Pakovich}

\begin{abstract}
With each holomorphic map $f:R\to\C\P^1$, where $R$
is a compact Riemann surface, one can associate a
combinatorial datum consisting of the genus $g$ of $R$,
the degree $n$ of $f$, the number $r$ of critical values
of $f$, and the $r$ partitions of $n$ given by the local
degrees of $f$ over these critical values. These quantities are related by the Riemann-Hurwitz formula, and the Hurwitz existence problem  asks whether
a combinatorial datum that fits this formula actually corresponds to some map $f$. In this paper, 
we use fiber products of holomorphic maps between
compact Riemann surfaces to establish general criteria
for the non-realizability of branch data. 
We also deduce from our results the theorem  of Halphen, proven in 1880,  concerning polynomial solutions of the equation 
$A(z)^a+B(z)^b=C(z)^c$, where $a,b,c$ are integers greater than one.

\end{abstract}

\maketitle

\section{Introduction}

Let $R$ be a compact Riemann surface of genus $g=g(R)$ and 
$f: R \rightarrow \C\P^1$ a holomorphic map of degree $n$. If  $z_1, z_2, ... ,z_r \in \C\P^1$
are the critical values of $f$, i.e. points $z\in 
\C\P^1$ for which $f^{-1}\{z\}$ contains fewer than  $n$ points, then for each $i,$ $1\leq i \leq r,$ the multiset
$\Pi_i=(\pi_{i,1}, \pi_{i,2},\, ... \, ,\pi_{i,p_i})$ of local degrees of $f$ at points of $f^{-1}\{z_i\}$ is a 
partition of $n$. Furthermore, it follows from the Riemann-Hurwitz formula that the equality 
\be \l{rh}  
\sum_{i=1}^r p_i=(r-2)n+2-2g
\ee 
holds. 
We call a collection of the form
$\Pi=(\Pi_1,\ldots,\Pi_r,n,g)$, where $n\geq1$, $g\geq0$,
and $r\geq0$ are integers, and $\Pi_i$, $1\leq i\leq r$,
are  partitions of $n$ satisfying \eqref{rh},
a \textit{branch datum}. All partitions occurring in a branch datum are assumed
to be non-trivial, meaning that they contain an entry
greater than one.
The Hurwitz existence problem is the following question:  
given a branch datum $\Pi = (\Pi_1, \ldots, \Pi_r, n, g)$, determine whether there exists a compact Riemann surface $R$ 
of genus $g$ and a holomorphic map $f: R \to \mathbb{C}\mathbb{P}^1$ such that $\Pi$ arises from $f$.  
If such a map exists, one says that $\Pi$ is \textit{realizable}; otherwise, that it is \textit{non-realizable}.

Note that the above problem is a special case of the broader problem of the existence of branched covering maps  
between closed connected  surfaces $\Sigma_1$ and $ \Sigma_2$, which goes back to Hurwitz (\cite{hur}). 
However, except when  $\Sigma_2$ is the sphere and $\Sigma_1$ is oriented, this  problem is either solved or can be reduced to this specific case (see \cite{eks}, \cite{pepe}), in which, by the Riemann existence theorem,  
  the existence of a branched covering map from $\Sigma_1$ to $\Sigma_2$  is equivalent to the existence of a holomorphic map $f: R \rightarrow \mathbb{C}\mathbb{P}^1$ between compact Riemann surfaces. Since our results are easier to formulate in terms of holomorphic maps, we use the corresponding formulation.

A large number of papers employing various techniques have been devoted to the Hurwitz problem (see \cite{adr}, \cite{b}- 
\hskip -0.1cm 
\cite{e}, \cite{fer}, \cite{ger}, \cite{hur}- 
\hskip -0.1cm \cite{kz}, \cite{ko}-
\hskip -0.1cm 
\cite{hu}, \cite{pap1}-
\hskip -0.1cm \cite{sx},  \cite{t}- 
\hskip -0.1cm \cite{zhu}), but the problem is still far from being solved. A comprehensive introduction to the topic can be found in \cite{pepe}.  
In this paper, we use fiber products of holomorphic maps between compact Riemann surfaces to prove several general results, providing a uniform explanation for the non-realizability of many known branch data and yielding numerous new examples. 

Our results are most conveniently formulated in terms of
orbifolds, and we begin by recalling the necessary definitions. {\it A Riemann surface orbifold} is a pair \linebreak $\f O=(R,\nu)$ consisting of a Riemann surface $R$ and a ramification function $\nu:R\rightarrow \mathbb N$ that takes the value $\nu(z)=1$ except at isolated points. 
For an orbifold $\f O=(R,\nu)$,  
 {\it the Euler characteristic} of $\f O$ is the number
$$ \chi(\f O)=\chi(R)+\sum_{z\in R}\left(\frac{1}{\nu(z)}-1\right),$$
{\it the set of  singular points} of $\f O$ is the  set 
$$c(\f O)=\{z_1,z_2, \dots, z_s, \dots \}=\{z\in R \mid \nu(z)>1\},$$ and  {\it the  signature} of $\f O$ is the multiset 
$$\nu(\f O)=\{\nu(z_1),\nu(z_2), \dots , \nu(z_s), \dots \}.$$ 
We say that an orbifold $\f O=(R,\nu)$  
is {\it non-ramified} if $\nu(z)\equiv 1.$ 

Let $\f O_1=(R_1,\nu_1)$  and $\f O_2=(R_2,\nu_2)$ be orbifolds
and let 
$f:\, R_1\rightarrow R_2$  be a holomorphic branched covering map. We say that $f:\,  \f O_1\rightarrow \f O_2$
is {\it a covering map} 
 between orbifolds
if for any $z\in R_1$ the equality 
\be \l{us} \nu_{2}(f(z))=\nu_{1}(z)\deg_zf\ee holds. 
For  a holomorphic map $f:R'\rightarrow R$ between compact Riemann surfaces and an  orbifold  $\f O=(R,\nu)$, we define its pullback  $f^*(\f O)=(R',\nu')$ by the formula 
\be \l{suu}  \nu'(z)= \frac{\nu(f(z))}{\gcd(\deg_zf, \nu(f(z)))},\ \ \ \ z\in R'.\ee
The map $f:f^*(\f O)\to\f O$ need not be a covering map
between orbifolds. However, if
$f:\f O_1\to\f O_2$ is a covering map between orbifolds,
then $f^*(\f O_2)=\f O_1$.

{\it A universal covering} of an orbifold ${\f O}$
is a covering map between orbifolds \linebreak  $\theta_{\f O}:\,
\tt {\f O}\rightarrow \f O$, where $\tt{\f  O}=(\tt R,\tt \nu)$ is a non-ramified orbifold  such that $\tt R$ is simply connected. 
If $\theta_{\f O}$ is such a map, then 
there exists a group $\Gamma_{\f O}$ of conformal automorphisms of $\tt R$ such that the equality 
$\theta_{\f O}(z_1)=\theta_{\f O}(z_2)$ holds for $z_1,z_2\in \tt R$ if and only if $z_1=\sigma(z_2)$ for some $\sigma\in \Gamma_{\f O}.$ 
A universal covering exists and is unique up to a conformal isomorphism of $\tt R$, unless $R = \C\P^1$ with a single ramified point, or $R = \C\P^1$ with two ramified points $z_1$, $z_2$ such that $\nu(z_1) \neq \nu(z_2)$.  
We say that an orbifold $\f O$ is \textit{good} if it has a universal covering, and \textit{bad} otherwise.

Furthermore, if $\f O = (R, \nu)$  is good and $R$ is compact, then the surface $\tt R$ is the unit disk $\mathbb D$ if and only if $\chi(\f O) < 0$, the complex plane $\C$ if and only if $\chi(\f O) = 0$, and the Riemann sphere $\C\P^1$ if and only if $\chi(\f O) > 0$. 
Finally, we recall that for a good orbifold $\f O$ on $\C\P^1$ the inequality 
$\chi(\f O)>0$ holds if and only if  $\nu(\f O)$ belongs to the  list 
 \be \l{lismain} \{d,d\}, \ \ d\geq 1,  \ \ \ \{2,2,d\}, \ \ d\geq 2,  \ \ \ \{2,3,3\}, \ \ \ \{2,3,4\}, \ \ \ \{2,3,5\}.\ee The corresponding universal coverings $\theta_{\f O}$ are well-known rational Galois coverings of $\C\P^1$ by $\C\P^1$ of degrees $d$, $2d$, $12$, $24$, and $60$, respectively, calculated for the first time by Klein (\cite{klein}). 
Most of the orbifolds considered in this paper are defined on $\mathbb{C}\mathbb{P}^1$. For this reason, abusing notation, we will usually
write simply ``an orbifold $\f O$'' when the orbifold
under consideration is defined on $\mathbb{C}\mathbb{P}^1$.

In the above notation, our main result is the following statement.

\begin{theorem}\label{main}
Let $f:\mathbb{CP}^1\to\mathbb{CP}^1$ be a rational
function, and let $\mathcal O$ be a good orbifold on
$\mathbb{CP}^1$.
Then the following statements hold.
\begin{enumerate}[label=\Roman*.,font=\normalfont\bfseries]
\item
The orbifold $f^*(\mathcal O)$ is good.

\item
If $\chi(f^*(\mathcal O))>0$, then $\chi(\mathcal O)>0$,
and, for suitably chosen universal coverings
$\theta_{\mathcal O}$ and $\theta_{f^*(\mathcal O)}$,
there exist rational functions $q$, $t$, and $w$
such that:
\begin{enumerate}[label=\arabic*),font=\normalfont]
\item
The equality
\[
f\circ\theta_{f^*(\mathcal O)}
=
\theta_{\mathcal O}\circ q
\]
holds.

\item
The equalities
\[
f=w\circ t,
\qquad
\theta_{\mathcal O}=w\circ\theta_{f^*(\mathcal O)}
\]
hold.

\item
The map $w:f^*(\mathcal O)\to\mathcal O$ is a covering
map between orbifolds. In particular,
\[
\deg w=\frac{\chi(f^*(\mathcal O))}{\chi(\mathcal O)}.
\]
\end{enumerate}
\end{enumerate}
\end{theorem}

Let us now mention some corollaries and particular cases of
Theorem~\ref{main}, along with examples.  

Assertion~I of the theorem offers a broad generalization of the following property of rational maps which follows easily from the existence of a universal covering:  
if  $f:\, \f O_1\rightarrow \f O_2$  is a covering map between orbifolds,  and $\f O_2$ is good, then 
$\f O_1$ is also good. 
The relationship between this property of covering maps and the Hurwitz problem, as well as its application in proving the non-realizability of certain branch data, was established in the paper \cite{pap1}. For example, 
 the branch datum 
 $$((2,2,2,2,1),(3,3,3),(3,3,3),9,0)$$ is non-realizable, since if a rational function $f$ with this branch datum existed, it would be a covering map $f:\, \f O_1\rightarrow \f O_2$, with $\f O_2$ good and $\f O_1$ bad. Here 
$\f O_2$ is defined by the equalities 
\be \l{bt} \nu_2(z_1)=2, \ \ \nu_2(z_2)=3, \ \ \nu_2(z_3)=3,\ee 
where $z_1,z_2,z_3$ are the critical values of $f$
corresponding, respectively, to the three partitions above, 
and $\f O_1$ is defined by the equality $\nu_1(z_0)=2$
at the unique non-critical point $z_0$ of $f$
in the fiber $f^{-1}(z_1)$.

The more general assertion I of Theorem \ref{main} allows us to construct many more examples of non-realizable data. As an example, let us consider the following well-known series of non-realizable branch data 
\be \l{iz0} \big((2^k),(2^k),(l,2k-l),2k,0\big), \ \ \ \ l\neq k\ee 
(hereinafter the symbol $u^v$ appearing in a partition means a string consisting of the number $u$ taken $v$ times). 
To see the non-realizability of \eqref{iz0} using
Theorem~\ref{main}, suppose that a rational function $f$
with this branch datum exists. Then, with $t=\max\{l,2k-l\}$, assigning weights $2,2,t$,
respectively, to the critical values corresponding to
the three partitions in \eqref{iz0} gives an orbifold
$\f O$ for which $f^*(\f O)$ has exactly one singular point.

Similarly, considering a suitable orbifold with signature
$\{2,3,5\}$, one can establish the non-realizability
of the branch data
\be \l{iz}
\big((2^l),(3^k),(5^m,s),n,0\big),
\qquad s\not\equiv0\pmod5,
\ee
previously proved in \cite{adr}, \cite{izm0}. 
More examples of non-realizable data obtained using
assertion~I of the theorem are given in Section~\ref{s31}.

The first interesting case to which assertion II of the theorem
applies is when $f^*(\f O)$ is non-ramified or equivalently when $\nu(f(z))$ divides $\deg_zf$ for all $z\in \C\P^1.$ 
In this case, $\theta_{f^*(\f O)}=id$ and both conditions 1) and 2) reduce to the existence of  a rational function $q$ such that the equality \be \l{crit} f = \theta_{\f O} \circ q\ee holds  (see Corollary \ref{t0}).  

It seems that the only previously known special case of this
statement is the following, corresponding to an orbifold
$\f O$ with signature $\{d,d\}$. 
Let $f:\C\P^1 \rightarrow \C\P^1$ be a rational function such that in the corresponding branch datum $\Pi$ all entries in $\Pi_1$ and $\Pi_2$ are divisible by an integer $d \geq 2$. Then, assuming that the critical values corresponding
to $\Pi_1$ and $\Pi_2$ are $0$ and $\infty$, respectively,
the equality $f=q^d$ holds for some rational function $q$. In this particular case, the statement follows immediately by
writing $f$ as a quotient of two polynomials, but even this elementary criterion
can  be used to establish the non-realizability of some
branch data, with the simplest example of this sort being $((2,2),(2,2),(1,3),4,0)$ (see \cite{pp}, Section 5, and \cite{wei} for more details and examples).

The more general Theorem \ref{main} provides a similar decomposability criterion \eqref{crit} for any finite group of rotations of the sphere. In particular, it   implies that a rational function such that $f^*(\f O)$ is non-ramified and $\deg f$ is not divisible by $\deg \theta_{\f O}$ cannot exist, 
allowing us to obtain large families of non-realizable branch data by using the necessary condition $ \deg \theta_{\f O} \mid\deg f$ alone. 
For example,  for the signature $\{2,3,3\}$, non-realizable branch data obtained in this way 
with the smallest possible values $n=18$ and $r=3$ are 
 $$((2^7,4),(3^6),(3^6),18,0),\ \ \ \ ((2^9),(3^4,6),(3^6),18,0),$$
while for the signature $\{2,3,5\}$ the first such data are 
$$\big((2^{43},4),(3^{30}),(5^{18}),90,0\big), \ \ \ \big((2^{45}),(3^{28},6),(5^{18}),90,0\big),$$ and $$\big((2^{45}),(3^{30}),(5^{16},10),90,0\big).$$  
For more examples of non-realizable branch data obtained from Theorem \ref{main} when $f^*(\f O)$ is non-ramified   we refer the reader to Section \ref{s32}.

 In the general case covered by assertion II, when $f^*(\f O)$
is ramified, Theorem~\ref{main} shows that the rational functions $f$ and
$\theta_{\f O}$ have a common compositional left factor
and determines its degree in terms of $\f O$ and  $f^*(\f O)$, where by a compositional left factor of a rational function  $f$ we mean any rational function $g$ such that  $f=g\circ h$ for some rational function $h.$

When $f^*(\f O)$ is ramified, assertion~II.1 of
Theorem~\ref{main} also allows us to establish
non-realizability by checking only a divisibility condition.
In particular, if $f$ is a rational function and
$\f O$ is a good orbifold such that $f^*(\f O)$ has
signature $\{a,a\}$, then Theorem~\ref{main} implies
that $a\cdot\deg f$ is divisible by $\deg\theta_{\f O}$.
This implies, for example, that the series of branch data
\be \l{such}
\big((2^{3l+3}),(3^{2l+2}),(5^{l+1},1,l),6(l+1),0\big),
\qquad l\equiv0\pmod2,
\ee
is non-realizable.
Indeed, if a rational function $f$ realizes \eqref{such},
assigning weights $2,3,5$, respectively, to the critical
values corresponding to the three partitions in
\eqref{such} gives an orbifold $\f O$ for which
$f^*(\f O)$ is either bad or has signature $\{5,5\}$. However, $5 \cdot 6(l+1)$ is not divisible by $\deg \theta_{\f O} = 60,$ in contradiction with Theorem \ref{main}.   

As a more subtle example of non-realizability provided by
Theorem~\ref{main}, we mention the following known series:
$$
\big((2^{k}),(2^{k-2},1,3),(k,k),2k,0\big),
\qquad k\geq2.
$$
In this case, the necessary divisibility condition is
satisfied; however, the equalities provided by
Theorem~\ref{main} are impossible
(see Section~\ref{s32}).

Note that the function $w$ may have degree one; this occurs when $\nu(f^*\f O) = \nu(\f O)$. In this case, the theorem yields no direct corollaries relevant to the Hurwitz problem. However, condition 1) remains meaningful and essentially describes a semiconjugacy relation between the rational functions $f$ and $q$ (see Section \ref{s30}).

Beyond its applications to the Hurwitz problem,
Theorem~\ref{main} has  an interesting 
consequence
concerning functional decompositions of rational functions.  
Let us recall that a rational function $f$ is called indecomposable if it cannot be represented as a composition of two rational functions of degree at least two. Otherwise,  $f$ is called decomposable. In group-theoretic terms, a rational function $f$ is indecomposable if and only if  the monodromy group of $f$ is primitive.  Decomposability is a rather delicate property that cannot
typically be seen from the branch datum of $f$ alone,
as the same branch datum can correspond to both primitive
and imprimitive monodromy groups.
Theorem~\ref{main} implies, however, that for certain
branch data, functions realizing those branch data are
necessarily decomposable, apart from explicitly described
exceptions (see Corollary~\ref{co}).

Another application of Theorem \ref{main} concerns plane graphs. Specifically, using geometric objects known as
``dessins d'enfants'' or ``constellations'', one can
interpret the realizability of a genus-zero branch datum
as the existence of a plane graph of a prescribed type
(see, e.g., \cite{ko}, \cite{lz}, \cite{hu}, \cite{pp}). Consequently, the results of this paper translate naturally
into graph-theoretic terms, yielding numerous results,
both old and new, on the existence of plane graphs. These applications are explored in a separate paper \cite{such}.

In short, our approach to proving Theorem~\ref{main}
is based on considering a component of the fiber product
of $f$ and $\theta_{\f O}$, which gives rise to the
commutative diagram
$$
\begin{CD}
R @>q>> \t R\\
@VV p V @VV \theta_{\f O} V\\
\C\P^1 @>f>> \C\P^1\,,
\end{CD}
$$
where $p$ and $q$ are holomorphic maps between 
Riemann surfaces.
An analysis of the ramification of the functions involved,
using Abhyankar's lemma, allows us to show that
$f^*(\f O)$ is good.
A further analysis of the relation between
the deck transformation groups of $p$ and $\theta_{\f O}$ enables us to prove the other statements of the theorem.

Note that in the examples of non-realizable branch data
obtained using assertion~II of Theorem~\ref{main}, the
obstacle to realizability lies in the ``forced
decomposability'' of the rational functions potentially
realizing these branch data. In particular, examples   obtained in this way, like other examples of  non-realizability  found so far, are consistent with the prime-degree conjecture proposed in \cite{eks},  which posits that any branch datum of prime degree  is realizable. 
Moreover, the examples derived from assertion I of
Theorem~\ref{main} cannot contradict this conjecture either,
since in each of these examples the degree $n$ has  a proper divisor among the weights of $\f O$.

Our results show that a condition which may indicate
the non-realizability of a branch datum $\Pi$ with $g=0$
is the existence of a good ramified orbifold $\f O$,
necessarily with positive Euler characteristic, such that
for a rational function $f$ that might realize this
datum, the pullback $f^*(\f O)$ has positive Euler
characteristic (and hence ``few'' singular points)
and a signature distinct from that of $\f O$.
Moreover, this condition is quite general. For example, among the 59 non-realizable
branch data with $g=0$, $r=3$, and $n\leq10$ listed
in  \cite{zh}, there are only seven for which
such an orbifold does not exist.

The findings in \cite{zh} demonstrate that a similar situation persists when passing from rational functions to holomorphic maps  from a torus to a sphere. Specifically, for all 30 non-realizable branch data
with $g=1$, $r=3$, and $n\leq20$ presented in \cite{zh},
the following holds: there exists a ramified orbifold
$\f O$ with $\chi(\f O)\geq0$ such that for a map $f$
that might realize this datum, the orbifold $f^*(\f O)$
is either non-ramified or has one or two singular points. 

The approach used to prove Theorem~\ref{main} is not directly
applicable to holomorphic maps from a torus to a sphere.
Nonetheless, a slight modification allows us to prove
the following result.

\bt \l{t3}  Let $R$ be a compact Riemann surface of genus one, let $f: R \to \C\P^1$ be a holomorphic map, and let $\f O$ be a good orbifold on $\C\P^1$ with $\chi(\f O) > 0$ such that $f^*(\f O)$ is non-ramified. Then there exist a rational function $w:\C\P^1 \rightarrow \C\P^1$ and a holomorphic map $t: R \to \C\P^1$ such that the equalities
\be 
f = w \circ t, \quad \theta_{\f O} = w \circ \theta_{w^*(\f O)} 
\ee
hold, and the signature of $w^*(\f O)$ is either $\{d, d\}$ for some $d \ge 1$, or $\{2, 2, 2\}$.

\et

Note that if in Theorem \ref{t3} the orbifold $\f O$ itself has signature $\{d, d\}$, $d \ge 1$, or $\{2, 2, 2\}$, then the conclusion holds trivially for $w = id$. However, if $\f O$ has a different signature, then Theorem \ref{t3} implies that $f$ is decomposable, and this fact can be used in the study of the Hurwitz problem. 
For example, using Theorem \ref{t3} one can establish  the non-realizability  of branch data   
$$\big((2^{k},k+3),(3^{k+1}),(3^{k+1}),3k+3,1\big), \ \ \ k\equiv 1\, ({\rm mod}\, 4),$$ 
and 
$$\big((2^{3k+6}),(3^{2k+4}),(3,9,6^{k}),6k+12,1\big), \ \ \ k\equiv 1\, ({\rm mod}\, 2),$$ 
  (see Section \ref{s42}).

This paper is organized as follows. In Section~\ref{s2},
we collect the necessary facts about orbifolds and
fiber products. In Section~\ref{s3}, we prove
Theorem~\ref{main} and discuss its applications to
the non-realizability of branch data and the
decomposability of rational functions.
In Section~\ref{s4}, we prove Theorem~\ref{t3}
and give some applications.

Finally, in Section~\ref{s5}, we use Theorem~\ref{main}
to deduce Halphen's theorem~\cite{hal} on polynomial solutions of the equation 
\[
A(z)^a+B(z)^b=C(z)^c,
\]
where $a$, $b$, and $c$ are integers greater than one.

\section{\l{s2} Orbifolds and fiber products} 
\subsection{Minimal holomorphic maps between orbifolds}
Throughout the paper, all holomorphic maps are assumed
to be non-constant.

In addition to the notion of a covering map between orbifolds given in the introduction, we will also use the notion of a holomorphic map between orbifolds. 
Let $\f O_1=(R_1,\nu_1)$  and $\f O_2=(R_2,\nu_2)$ be orbifolds
and let 
$f:\, R_1\rightarrow R_2$  be a holomorphic branched covering map. We say that $f:\, \f O_1\rightarrow \f O_2$ is a holomorphic map
between orbifolds if, for every $z\in R_1$, the weaker condition
\be \l{uuss} \nu_2(f(z))\mid \nu_1(z)\deg_zf\ee
holds in place of equality \eqref{us}.

If $\f O_1$ and $\f O_2$ are good, with universal coverings
$\theta_i:\t R_i\to R_i$, $i=1,2$, then a covering map
$f:\f O_1\to\f O_2$ lifts to a conformal isomorphism
$\t f$ making the diagram
\be  \l{dsa} 
\begin{CD}
\t R_1 @>\t f>> \t R_2\\
@VV  \theta_1 V @VV  \theta_2 V\\ 
R_1 @>f >> R_2 
\end{CD}
\ee
commutative. On the other hand, if $f$ satisfies \eqref{uuss}, it still
admits a holomorphic lift $\t f$ making diagram \eqref{dsa}
commutative, but this lift need not be an isomorphism
(see \cite{semi}, Proposition 3.1). This explains the term ``holomorphic map'' between orbifolds.

If $f:\,  \f O_1\rightarrow \f O_2$ is a covering map between orbifolds with compact $R_1$ and $R_2$, then  the Riemann-Hurwitz 
formula implies that 
\be \l{rhor} \chi(\f O_1)=d \chi(\f O_2), \ee
where $d=\deg f$. 
For holomorphic maps the following statement is true (see \cite{semi}, Proposition 3.2). 

\bp \l{p1} Let $f:\, \f O_1\rightarrow \f O_2$ be a holomorphic map between orbifolds with compact $R_1$ and $R_2$.
Then 
\be \l{iioopp} \chi(\f O_1)\leq \chi(\f O_2)\,\deg f, \ee and equality 
holds if and only if $f:\, \f O_1\rightarrow \f O_2$ is a covering map between orbifolds. 
\ep

Let $R_1$, $R_2$ be Riemann surfaces and 
$f:\, R_1\rightarrow R_2$ a holomorphic branched covering map. Assume that $R_2$ is equipped with a ramification function
$\nu_2$. To make $f$ a holomorphic map between the orbifolds
$\f O_1=(R_1,\nu_1)$ and $\f O_2=(R_2,\nu_2)$, we must
choose $\nu_1$ to satisfy condition \eqref{uuss}.
It is easy to see that, for every $z\in R_1$, the smallest
possible value of $\nu_1(z)$ is determined by the equality 
\be \l{rys} \nu_{2}(f(z))=\nu_{1}(z)\GCD(\deg_zf, \nu_{2}(f(z))).\ee 
If \eqref{rys} holds for every $z\in R_1$, we
say that $f$ is  a  minimal holomorphic  map
between orbifolds 
$\f O_1=(R_1,\nu_1)$ and $\f O_2=(R_2,\nu_2)$.

It follows from the definition that for any orbifold $\f O = (R, \nu)$ and any holomorphic branched covering map $f\colon R' \to R$, there exists a unique orbifold structure $\nu'$ on $R'$ such that $f$ becomes a minimal holomorphic map between the orbifolds $\f O' = (R', \nu')$ and $\f O = (R, \nu)$. We denote the corresponding orbifold $\f O'$ by $f^*\f O$, consistent with the notation and formula \eqref{suu} from the introduction.
Let us emphasize that for any minimal holomorphic map $f\colon \f O_1 \to \f O_2$, the equality \be \l{zzxx} \f O_1 = f^*(\f O_2) \ee holds simply by definition.
Note also that any covering map between orbifolds $f\colon \f O_1 \to \f O_2$ is a minimal holomorphic map. In particular, equality \eqref{zzxx} holds.

\vskip 0.2cm
Minimal holomorphic maps between orbifolds possess the following fundamental property (see \cite{semi}, Theorem 4.1).

\bt \l{serrr} Let $f:\, R^{\prime\prime} \rightarrow R^{\prime}$ and $g:\, R^{\prime} \rightarrow R$ be holomorphic branched covering maps, and  $\f O=(R,\nu)$ an orbifold. 
Then 
$$(g\circ f)^*(\f O)= f^*(g^*(\f O)).\eqno{\Box}$$
\et

Theorem \ref{serrr} implies in particular the following corollaries (see   \cite{semi}, Coro- \linebreak llary 4.1 and Corollary 4.2).

\bc \l{serka0} Let $f:\, \f O_1\rightarrow \f O^{\prime}$ and $g:\, \f O^{\prime}\rightarrow \f O_2$ be minimal holomorphic maps (resp. covering maps) between orbifolds.
Then  $g\circ f:\, \f O_1\rightarrow \f O_2$ is  a minimal holomorphic map (resp. covering map). \qed
\ec

\bc \l{indu2}  Let $f:\, R_1 \rightarrow R^{\prime}$ and $g:\, R^{\prime} \rightarrow R_2$ be holomorphic branched covering maps, and  $\f O_1=(R_1,\nu_1)$ and  $\f O_2=(R_2,\nu_2)$
orbifolds. Assume that \linebreak $g\circ f:\, \f O_1\rightarrow \f O_2$ is  a minimal holomorphic map between orbifolds (resp. a co\-vering map). Then  $g:\, g^*(\f O_2)\rightarrow \f O_2$  and  $f:\, \f O_1\rightarrow g^*(\f O_2) $ are minimal holomorphic maps (resp. covering maps). \qed
\ec

For   orbifolds $\f O_1=(R_1,\nu_1)$  and $\f O_2=(R_2,\nu_2)$, we   
write \be \l{elki} \f O_1\preceq \f O_2 \ee 
if $R_1=R_2$, and for any $z\in R_1$ the condition $$\nu_1(z)\mid \nu_2(z)$$ holds.

For a Riemann surface $R$, we denote by $\f R$ the
non-ramified orbifold on $R$. When $R=\C\P^1$,
we simply retain the notation $\C\P^1$ for the
orbifold $\f R$.

\subsection{Fiber products} 
Let $f:\, C\rightarrow \C\P^1$ and $g:\, D\rightarrow \C\P^1$ be holomorphic maps between compact Riemann surfaces. 
The collection
\be \l{nota}
(C,f)\times(D,g)
=\bigl\{\{R_j,p_j,q_j\}:1\leq j\leq n(f,g)\bigr\},
\ee
where $n(f,g)$ is a positive  integer  and $R_j$ are compact Riemann surfaces provided with holomorphic maps
$$p_j:\, R_j\rightarrow C, \ \ \ q_j:\, R_j\rightarrow D, \ \ \ 1\leq j \leq n(f,g),$$
is called the {\it fiber product} of  $f$ and $g$ if \be \l{pes} f\circ p_j=g\circ q_j, \ \ \ 1\leq j \leq n(f,g),\ee 
and for any holomorphic maps $p:\, R\rightarrow C,$  $q:\, R\rightarrow D$
between compact Riemann surfaces satisfying 
\be \l{mm} f\circ p=g\circ q\ee there exist a uniquely defined  index $j$ and 
a holomorphic map $w:\, R\rightarrow R_j$ such that
\be \l{univ} p= p_j\circ  w, \ \ \ q= q_j\circ w.\ee
The fiber product exists and is unique up to natural
isomorphisms. Here an isomorphism between components
$\{R,p,q\}$ and $\{R',p',q'\}$ is a biholomorphic map
$\alpha:R\to R'$ such that
$p=p'\circ\alpha$ and $q=q'\circ\alpha$.

The fiber product can be described by the following
algebro-geometric construction. Let us consider the algebraic curve
\be \l{ccuurr}
L=\{(x,y)\in C\times D\,\vert\,f(x)=g(y)\}.
\ee
Let us denote by $L_j$, $1\leq j\leq n(f,g)$, the irreducible
components of $L$, by $R_j$ their desingularizations, and by
$$
\pi_j:\,R_j\rightarrow L_j,\qquad 1\leq j\leq n(f,g),
$$
the desingularization maps.
Let $x$ and $y$ denote the projections of $C\times D$
onto its factors. Then the compositions
$$
p_j=x\circ\pi_j:\,R_j\rightarrow C,\qquad
q_j=y\circ\pi_j:\,R_j\rightarrow D,\qquad
1\leq j\leq n(f,g),
$$
are holomorphic maps, and the collection
$\bigl\{\{R_j,p_j,q_j\}:1\leq j\leq n(f,g)\bigr\}$
is the fiber product of $f$ and $g$. 
Abusing notation, we refer to the Riemann surfaces $R_j$,
$1\leq j\leq n(f,g)$, as the irreducible components of the
fiber product of $f$ and $g$, and say that the fiber product
is irreducible if $n(f,g)=1$.

It follows from the definition  that for every $j,$ $1\leq j \leq n(f,g),$ the functions $p_j,q_j$  have no {\it non-trivial common compositional  right factor} in the following sense: 
the equalities 
$$ p_j= \tt p\circ  t, \ \ \ q_j= \tt q\circ t,$$ where $$t:\, {R}_j \rightarrow \t{R} , \ \ \ \tt p:\, \t{R} \rightarrow C, \ \ \ \tt q:\, \t{R}  \rightarrow D$$ are holomorphic maps between compact Riemann surfaces, imply that $\deg t=1.$

Denoting by $\f M(R)$ the field of meromorphic functions
on a Riemann surface $R$, we can restate this condition
as the equality
$$
p_j^*\f M(C)\cdot q_j^*\f M(D)=\f M(R_j),
$$
where the product denotes the compositum 
of these two subfields. Indeed, by the correspondence between
holomorphic maps and extensions of meromorphic function
fields, a non-trivial common compositional right factor
exists precisely when both subfields are contained in
a proper subfield of $\f M(R_j)$.

Conversely, if $p$ and $q$ satisfy \eqref{mm} and have no
non-trivial common compositional right factor, then
$$
p=p_j\circ t,\qquad q=q_j\circ t
$$
for some $j$, $1\leq j\leq n(f,g)$, and an isomorphism
$t:\,R\rightarrow R_j$.

Since the pairs $(p_j,q_j)$,
$1\leq j\leq n(f,g)$, parametrize the components of
the curve \eqref{ccuurr}, we have: 
\be \l{ii}
\sum_j\deg p_j=\deg g,\qquad
\sum_j\deg q_j=\deg f. 
\ee

\subsection{Functional equations and orbifolds}

 With each holomorphic map \linebreak $h:\, R_1\rightarrow R_2$ between compact Riemann surfaces, 
one can associate two orbifolds $\f O_1^h=(R_1,\nu_1^h)$ and 
$\f O_2^h=(R_2,\nu_2^h)$ in a natural way, setting $\nu_2^h(z)$  
equal to the least common multiple of local degrees of $h$ at the points 
of the preimage $h^{-1}\{z\}$, and $$\nu_1^h(z)=\frac{\nu_2^h(h(z))}{\deg_zh}.$$ By construction, 
 $$h:\, \f O_1^h\rightarrow \f O_2^h$$ 
is a covering map between orbifolds.
It is easy to see that these orbifolds are minimal
in the following sense. For any covering map   $h:\, \f O_1\rightarrow \f O_2$, we have:
\be \l{elki+} \f O_1^h\preceq \f O_1, \ \ \ \f O_2^h\preceq \f O_2.\ee Note that the orbifolds $\f O_1^h$ and $ \f O_2^h$ 
 are good  (see \cite{semi}, Lemma 4.2). 
Note also that for the universal covering $\theta_{\f O}$ of an orbifold $\f O$ of positive Euler characteristic the equalities 
$$\f O_2^{\theta_{\f O}}=\f O, \ \ \ \  \f O_1^{\theta_{\f O}}=\C\P^1$$ hold.

The orbifolds defined above are useful for the study of the functional equation 
\eqref{mm}, where  
$$p:\, R\rightarrow C,\ \ f:\, C\rightarrow \C\P^1,\  \ q:\, R\rightarrow D,\  \ g:\, D\rightarrow \C\P^1$$ 
are holomorphic maps between compact Riemann surfaces.
We will usually write this equation in the form of
a commutative diagram
\be \l{m}
\begin{CD}
R @>q>> D\\
@VV {p} V @VV g V\\ 
C @>f >> \C\P^1\,.
\end{CD}
\ee
Specifically, the following statement holds (see \cite{semi}, Theorem 4.2).

\bt \l{goodt} Let $f,p,g,q$ be holomorphic maps between compact Riemann surfaces such that the  diagram \eqref{m} commutes,  the fiber product of $f$ and $g$ has a unique component, 
and $p$ and $q$ have no non-trivial common compositional right factor. 
Then the commutative diagram 
\be 
\begin{CD}
\f O_1^p @>q>> \f O_1^g\\
@VV p V @VV g V\\ 
\f O_2^p @>f >> \f O_2^g\ 
\end{CD}
\ee
consists of minimal holomorphic  maps between orbifolds. \qed  
\et

Since the vertical arrows in the above diagram are covering
maps between orbifolds, they are automatically minimal
holomorphic maps. The nontrivial content of the theorem
is expressed by the equalities 
$$\f O_2^p=f^*( \f O_2^g), \ \ \ \ \  \f O_1^p=q^*( \f O_1^g),$$
which describe the orbifolds $\f O_2^p$ and $\f O_1^p$ as pullbacks.

In this paper, we will use a version of Theorem~\ref{goodt}
that also applies to reducible fiber products, under
an additional assumption on $g$.

\section{\l{s3} Proof and applications of Theorem~\ref{main}}
\subsection{\l{s30} Proof of Theorem~\ref{main}}
Let $p:\, R\rightarrow C$ be a holomorphic map between compact Riemann surfaces. 
Let us recall that $p$ is called a \textit{Galois covering} if its automorphism group
$$\Aut(R,p)=\{\sigma\in \Aut(R)\,:\, p\circ\sigma=p\}$$
acts transitively on fibers of $p$.  
Equivalently, $p$ is a Galois covering if the field extension 
$ \f M(R)/p^*\f M(C)$ is a Galois extension. 
If $p$ is a Galois covering, then
$$
\Aut(R,p)\cong
{\rm Gal}\left(\f M(R)/p^*\f M(C)\right).
$$
Since $\Aut(R,p)$ acts freely on every fiber over a
non-critical value of $p$, the map $p$ is a Galois covering
if and only if
\be \l{dega}
|\Aut(R,p)|=\deg p.
\ee

Let $h:R\to C$ be a holomorphic map between compact Riemann surfaces.
We say that $h$ is {\it uniform} if $\f O_1^h$ is non-ramified.
Equivalently, the local degree of $h$ is constant on each fiber. Note that every Galois covering is uniform, but the converse
is not true in general.

We start by proving two results of independent interest.
The first is a version of Theorem~\ref{goodt} that applies
to possibly reducible fiber products under the additional
assumption that $g$ is uniform. The proof of the second
extends the group-theoretic argument used in the proof
of Theorem~5.1 in \cite{semi}.

\bt \l{t5}
Let $f,p,g,q$ be holomorphic maps between compact Riemann
surfaces such that diagram \eqref{m} commutes, and $p$
and $q$ have no non-trivial common compositional right
factor. Assume that $g$ is uniform. Then $p$ is uniform,
and the commutative diagram
$$
\begin{CD}
\f O_1^p @>q>> \f O_1^g\\
@VV p V @VV g V\\
\f O_2^p @>f>> \f O_2^g
\end{CD}
$$
consists of minimal holomorphic maps between orbifolds.
In particular, if $f^*(\f O_2^g)$ is non-ramified,
then $p$ is unbranched.
\et

\pr Let us set
$$
F=f\circ p=g\circ q.
$$
Since $p$ and $q$ have no non-trivial common compositional
right factor, by Abhyankar's lemma
(see e.g. \cite{sti}, Theorem~3.9.1), for every point $t_0$
of $R$ we have
\be \l{abj}
e_F(t_0)=
{\rm lcm}\big(e_f(p(t_0)),e_g(q(t_0))\big),
\ee
where, for brevity, we temporarily write $e_h(t)=\deg_t h$
for a holomorphic map $h$ and a point $t$ in its domain. 
Since $g$ is uniform,
\be \l{imb}
e_g(q(t_0))
=\nu_2^g\big(g(q(t_0))\big)
=\nu_2^g\big(f(p(t_0))\big).
\ee
It follows now from \eqref{abj}, \eqref{imb}, and the chain
rule that
\be 
e_p(t_0)
=\frac{e_F(t_0)}{e_f(p(t_0))}=\frac{{\rm lcm}\big(e_f(p(t_0)),e_g(q(t_0))\big)}
        {e_f(p(t_0))}
=\frac{\nu_2^g\big(f(p(t_0))\big)}
        {\gcd\big(e_f(p(t_0)),\nu_2^g(f(p(t_0)))\big)}.
\ee
The right-hand side is the value at $p(t_0)$ of the
ramification function of $f^*(\f O_2^g)$. Therefore,
$
p:\f R\longrightarrow f^*(\f O_2^g)
$
is a covering map between orbifolds. Since its source is non-ramified, it follows from
\eqref{elki+} that $p$ is uniform. 
Furthermore, since \linebreak $p:\f R\to\f O_2^p$ is also
a covering map, \eqref{us} yields
\be \l{ash} \f O_2^p=f^*(\f O_2^g).\ee
In particular, if $f^*(\f O_2^g)$ is non-ramified,
then $p$ is unbranched.
Equality \eqref{ash} also shows that
$f:\f O_2^p\to\f O_2^g$ is minimal by the definition
of the pullback orbifold.

 Finally, since both $\f O_1^p$
and $\f O_1^g$ are non-ramified, the map $q:\f O_1^p\rightarrow \f O_1^g$ is trivially minimal, while the vertical
maps are covering maps by construction. 
\qed

\bt \l{gc}
Let $f,p,g,q$ be holomorphic maps between compact Riemann
surfaces such that diagram \eqref{m} commutes, and $p$
and $q$ have no non-trivial common compositional right
factor. Assume that $g$ is a Galois covering. Then there
exists an injective homomorphism
$$
\phi:\Aut(R,p)\longrightarrow\Aut(D,g)
$$
such that
\be \l{hh}
q\circ\sigma=\phi(\sigma)\circ q,
\qquad \sigma\in\Aut(R,p).
\ee
\et

\pr 
Clearly, the commutativity of 
\eqref{m} implies that
for every $\sigma\in\Aut(R,p)$ the equality
\be\l{as}
g\circ(q\circ\sigma)=g\circ q
\ee
holds. On the other hand, for the fiber product of $g$
with itself, the maps $p_j,q_j$ in \eqref{nota}
can be taken as
$$
p_j=\mu_j,\qquad q_j=id,\qquad
\mu_j\in\Aut(D,g).
$$
Indeed, $p_j,q_j$ defined in this way satisfy \eqref{pes}
and have no non-trivial common compositional right factor.
Moreover, since
$$
\sum_j\deg\mu_j=|\Aut(D,g)|=\deg g
$$
by \eqref{dega}, these maps exhaust all components
in \eqref{nota} by \eqref{ii}.
Thus, the universality property
$$
q\circ\sigma=p_j\circ w,\qquad q=q_j\circ w
$$
reduces to equality \eqref{hh} for some
$\phi(\sigma)\in\Aut(D,g)$.
Finally, one can easily see that the correspondence
$\sigma\mapsto\phi(\sigma)$ defines a homomorphism of groups.

To complete the proof, we only need to show that the group
$\Ker\phi$ is trivial. 
To this end, observe that if $\Ker\phi$ is non-trivial,
then the set of meromorphic functions $h$ on $R$ satisfying
the condition
$$
h\circ\sigma=h,\qquad \sigma\in\Ker\phi,
$$
forms a subfield $k$ of $\f M(R)$ distinct from $\f M(R)$.
By \eqref{hh}, we have $q^*\f M(D)\subseteq k$. 
Moreover, $p^*\f M(C)\subseteq k$, since $\Ker\phi$ 
is a subgroup of $\Aut(R,p)$. Thus, if $\Ker\phi$ is non-trivial,
$$
p^*\f M(C)\cdot q^*\f M(D)\neq\f M(R),
$$
in contradiction with the assumption that $p$ and $q$
have no non-trivial common compositional right factor. \qed

The following lemma is the final auxiliary result needed to prove Theorem~\ref{main}. It follows directly from the definition of a covering map and two well-known facts: an orbifold on the sphere with positive Euler characteristic has at most three singular points, and any ordered triple of distinct points on the sphere can be mapped to any other such triple by a M\"obius transformation.

\bl \l{lmo}
Let $\f O_1$ and $\f O_2$ be orbifolds on $\C\P^1$
with positive Euler characteristic. Then
$\nu(\f O_1)=\nu(\f O_2)$ if and only if there exists
a M\"obius transformation $w$ such that
$w:\f O_1\to\f O_2$ is a covering map between orbifolds. \qed 
\el

\vskip 0.2cm
\noindent{\it Proof of Theorem~\ref{main}.}
If $\chi(f^*(\f O))\leq0$, then $f^*(\f O)$
is good, since every bad orbifold on the sphere has
positive Euler characteristic. It therefore remains
to consider the case $\chi(f^*(\f O))>0$.
Proposition~\ref{p1} gives
\[
0<\chi(f^*(\f O))\leq\chi(\f O)\,\deg f.
\]
Thus, $\chi(\f O)>0$ and $\theta_{\f O}$ is
a rational Galois covering.

Choose a component $\{R,p,q\}$ of the fiber product
of $f$ and $g=\theta_{\f O}$.
Since $\theta_{\f O}$ is a universal covering of $\f O$,
it is uniform and $\f O_2^{\theta_{\f O}}=\f O$.
Applying Theorem~\ref{t5},
we conclude that $p$ is uniform, so that $\f O_1^p=\f R$, and
\[
\f O_2^p=f^*(\f O).
\]
Thus, the canonical covering map
$p:\f O_1^p\to\f O_2^p$ takes the form
\be \l{u}
p:\f R\to f^*(\f O).
\ee Moreover, by \eqref{rhor},  
\[
\chi(\f R)=2-2g(R)=\chi(f^*(\f O))\,\deg p>0.
\]
Hence $R$ is a sphere and therefore simply connected.
Since $\f R$ is non-ramified, we see that the map \eqref{u} is
a universal covering of $f^*(\f O)$.
Therefore, the orbifold $f^*(\f O)$ is good, and $p$ is
a Galois covering. 
Moreover, identifying $R$ with $\C\P^1$
and $p$ with $\theta_{f^*(\f O)}$, we also obtain assertion~II.1:
\be\l{eqr}
f\circ\theta_{f^*(\f O)}=\theta_{\f O}\circ q.
\ee 

To prove the remaining assertions, set
\[
H=\Aut(\C\P^1,\theta_{f^*(\f O)}),\qquad G=\Aut(\C\P^1,\theta_{\f O}), 
\] 
and let $$K=\phi(H)\subseteq G,$$ where $\phi$ is the embedding
provided by Theorem~\ref{gc}. Identifying the quotient $\C\P^1/K$ with $\C\P^1$,
let $v:\C\P^1\to\C\P^1$ be the quotient map.
Since $K\subseteq G$, there exists a rational function $h$
such that
$
\theta_{\f O}=h\circ v.
$
Furthermore,  setting $\f O'=h^*(\f O)$, we see that both $v:\C\P^1\to\f O'$ and
$h:\f O'\to\f O$ are covering maps between orbifolds by Corollary~\ref{indu2}.  
In particular, $v$ is a universal covering of $\f O'$.

Since $v$ 
is the quotient map by $K$,  the map $v\circ q$ is $H$-invariant by \eqref{hh}. Since $\theta_{f^*(\f O)}$
is the quotient map by $H$, this implies that there exists a rational map
$t_0:\C\P^1\to \C\P^1$ such that
\[
v\circ q=t_0\circ \theta_{f^*(\f O)}.
\]
It follows that $$f\circ \theta_{f^*(\f O)}=\theta_{\f O}\circ q=h\circ v\circ q=h\circ t_0\circ \theta_{f^*(\f O)},$$ and hence
$f=h\circ t_0$.

Since $\theta_{f^*(\f O)}:\C\P^1\to f^*(\f O)$ and
$v:\C\P^1\to\f O'$ are universal coverings with isomorphic
deck transformation groups $H$ and $K$, respectively,
the classification of finite groups of M\"obius
transformations implies that
\[
\nu(f^*(\f O))=\nu(\f O').
\] 
Lemma~\ref{lmo} therefore provides a M\"obius transformation
$\mu:f^*(\f O)\to\f O'$ which is an isomorphism of orbifolds.
Set $$w=h\circ\mu, \qquad t=\mu^{-1}\circ t_0.$$
Then $f=w\circ t$, and \be \l{coma} w:f^*(\f O)\to\f O\ee is a covering map
by Corollary~\ref{serka0}.

Finally, since $\theta_{f^*(\f O)}$ and $\mu^{-1}\circ v$ are universal coverings of
$f^*(\f O)$, there is a M\"obius transformation $\alpha$ such that
$$\theta_{f^*(\f O)}\circ\alpha=\mu^{-1}\circ v.$$ Replacing now 
$\theta_{f^*(\f O)}$ by $\theta_{f^*(\f O)}\circ\alpha$ and $q$ by $q\circ\alpha$
preserves \eqref{eqr} and gives 
\be\l{cond}
f=w\circ t,\qquad
\theta_{\f O}=h\circ v=w\circ\mu^{-1}\circ v=w\circ\theta_{f^*(\f O)}.
\ee
Moreover, 
\[
\chi(f^*(\f O))=\chi(\f O)\,\deg w,
\]
by \eqref{rhor}, 
which completes the proof.
\qed

\bc \l{t0} Let $f:\C\P^1 \rightarrow \C\P^1$ be a rational function, and let $\f O$ be a good orbifold such that $f^*(\f O)$ is non-ramified. Then $\chi(\f O) > 0$, and $f = \theta_{\f O} \circ q$ for some rational function $q$. In particular, $\deg f$ is divisible by $\deg \theta_{\f O}$.
\ec
\pr
Since $f^*(\f O)$ is non-ramified, we have
$\chi(f^*(\f O))=2>0$ and may take
$\theta_{f^*(\f O)}=id$. Then 
assertion~II.1 of Theorem~\ref{main} gives
$\chi(\f O)>0$ and $f=\theta_{\f O}\circ q$.
\qed
\vskip 0.2cm
Note that if
\be \l{noth} 
\nu(f^*(\f O))=\nu(\f O),
\ee 
then the function $w$ in Theorem~\ref{main}
has degree one, and the map \eqref{coma}
is an isomorphism of orbifolds.
Since $f:f^*(\f O)\to\f O$ is a minimal holomorphic map by Theorem \ref{t5}, it follows
by Corollary~\ref{serka0} that the rational function
$\h f=f\circ w^{-1}$ is a minimal holomorphic self-map
$\h f:\f O\to\f O$.

For $\deg f\geq2$ and ramified $\f O$, such rational
functions are called \textit{generalized Latt\`es maps},
by analogy with ordinary Latt\`es maps, for which the
self-map is required to be a covering map.
Generalized Latt\`es maps play a crucial role in the study
of semiconjugate rational functions
(see \cite{semi}, \cite{lattes}).
In particular, since \eqref{cond} gives
$$w^{-1}\circ\theta_{\f O}=\theta_{f^*(\f O)},$$
equation \eqref{eqr} can be written as
\[
\h f\circ\theta_{\f O}=\theta_{\f O}\circ q.
\]
It would be interesting to understand which branch data
potentially corresponding to generalized Latt\`es maps
are realizable. For ordinary Latt\`es maps, this question
was resolved in \cite{pap1}.

\subsection{\l{s32} Non-realizability from factorizations}
 We begin with applications of  Corollary~\ref{t0},  which gives both a divisibility
condition on the degree and a prescribed compositional
left factor.

For integers $a,b,c>1$, we set  $\chi(a,b,c)=\chi(\f O)$, where $\f O$ is any orbifold on $\C\P^1$ with the signature
$\{a,b,c\}$. When $\chi(a,b,c)>0$, we also define
$$
n(a,b,c)=\deg\theta_{\f O},
\qquad
l(a,b,c)=\mathrm{lcm}(a,b,c).
$$
Since the Euler characteristic of a non-ramified orbifold on $\C\P^1$ equals two, we have
$$
\chi(a,b,c) = \frac{2}{n(a,b,c)},
$$
by \eqref{rhor}.  In addition, it is easy to see that
\be \l{by}
l(a,b,c)=\frac{n(a,b,c)}{2},
\ee
unless $(a,b,c)$ is a permutation of $(2,2,m)$ with
$m$ odd.

For a partition $\Pi=(\pi_1,\pi_2,\ldots,\pi_\ell)$
of a number $n$, we define the number $d(\Pi)$ by
the formula
\be \l{dp}
d(\Pi)=n-\ell.
\ee 
In this notation, the Riemann--Hurwitz condition
\eqref{rh} for a genus-zero branch datum
$\Pi=(\Pi_1,\Pi_2,\dots,\Pi_r,n,0)$ with $r\geq3$ becomes
\be \l{rh3}
p_1+p_2+p_3=n+2+\sum_{i=4}^r d(\Pi_i),
\ee
where $p_i$ is the number of entries of $\Pi_i$.

\bc\l{c1} 
Let $\Pi=(\Pi_1,\Pi_2,\dots, \Pi_r,n,0)$ be a branch datum. Assume that there exists 
 a triple of integers greater than one $(a,b,c)$
 such that all entries of $\Pi_1$ are divisible by $a$,  all entries of $\Pi_2$ are divisible by $b$, and all entries of $\Pi_3$ are divisible by $c$. Then $\chi(a,b,c)>0$. Moreover, $\Pi$ is non-realizable
unless $n$ is divisible by $n(a,b,c)$.
\ec
\pr 
By \eqref{rh3} and the divisibility assumptions,
$$
n+2\leq p_1+p_2+p_3
\leq \frac{n}{a}+\frac{n}{b}+\frac{n}{c},
$$
so
$$
\chi(a,b,c)
=\frac1a+\frac1b+\frac1c-1
\geq\frac2n>0.
$$

Furthermore, suppose that $f$ realizes $\Pi$, and let
$\f O$ be the orbifold defined by
$$
\nu(z_1)=a,\qquad \nu(z_2)=b,\qquad \nu(z_3)=c,
$$
where $z_1,z_2,z_3$ are the critical values of $f$
corresponding, respectively, to the partitions
$\Pi_1,\Pi_2,\Pi_3$. 
Then $f^*(\f O)$ is non-ramified, so $n(a,b,c)$ divides
$\deg f=n$ by Corollary~\ref{t0}.
\qed

\vskip 0.2cm

The following proposition demonstrates the existence of branch data whose non-realizability results from Corollary \ref{c1} and provides a method for the practical construction of such data.

\bp  \l{proo} 
Let $(a,b,c)$ be a triple of integers greater than one
satisfying $\chi(a,b,c)>0$ that is not a permutation of
$(2,2,m)$ for an odd integer $m$. Then for any integer of the form \be \l{w1} n=\frac{n(a,b,c)}{2}k,\ \ \ k\geq 3,\ee where  $k$ is odd, there exists  a non-realizable branch datum $\Pi=(\Pi_1,\Pi_2,\dots, \Pi_r,n,0)$  such that  all entries of $\Pi_1$ are divisible by $a$,  all entries of $\Pi_2$ are divisible by $b$, and  all entries of $\Pi_3$ are divisible by $c$. 
\ep

\begin{proof} Let $n$ be given by \eqref{w1}, and let $r\geq3$
be an integer. If $r>3$, fix non-trivial partitions
$\Pi_4,\ldots,\Pi_r$ of $n$. 
For positive integers $u_i$, $v_j$, and $w_e$, set
\be \l{0}
\Pi_1=(au_1,\ldots,au_{p_1}),\qquad
\Pi_2=(bv_1,\ldots,bv_{p_2}),\qquad
\Pi_3=(cw_1,\ldots,cw_{p_3}).
\ee
We first observe that the collection
$\Pi=(\Pi_1,\Pi_2,\ldots,\Pi_r,n,0)$ is a branch datum
if and only if these integers satisfy the system
\be \l{1}
\sum_{i=1}^{p_1}u_i=\frac{n}{a},\qquad
\sum_{j=1}^{p_2}v_j=\frac{n}{b},\qquad
\sum_{e=1}^{p_3}w_e=\frac{n}{c},
\ee
\be \l{22}
\sum_{i=1}^{p_1}(u_i-1)
+\sum_{j=1}^{p_2}(v_j-1)
+\sum_{e=1}^{p_3}(w_e-1)
=k-2-\sum_{i=4}^r d(\Pi_i).
\ee

Indeed, equalities \eqref{1} mean that
$\Pi_1,\Pi_2,\Pi_3$ are partitions of $n$.
Under these equalities, we have
\be \l{imp}
\begin{split}
\sum_{i=1}^{p_1}(u_i-1)
+\sum_{j=1}^{p_2}(v_j-1)
+\sum_{e=1}^{p_3}(w_e-1)
&=\frac{n}{a}+\frac{n}{b}+\frac{n}{c}
  -(p_1+p_2+p_3)\\
&=n+k-(p_1+p_2+p_3),
\end{split}
\ee
since $n\chi(a,b,c)=k$ by \eqref{w1}.
Thus, \eqref{22} is equivalent to \eqref{rh3}.

Since $k$ is odd, all branch data satisfying
\eqref{1} and \eqref{22} are non-realizable by
Corollary~\ref{c1}.
Finally, it is easy to see that if $r=3$, or, more
generally, if $r\geq3$ and $\Pi_4,\ldots,\Pi_r$ are
arbitrary non-trivial partitions of $n$ given by
\eqref{w1} satisfying
\be \l{w2}
\sum_{i=4}^r d(\Pi_i)\leq k-2,
\ee
then the system \eqref{1}, \eqref{22} has solutions
in positive integers $u_i$, $v_j$, $w_e$. Indeed, condition \eqref{w2} implies that the right-hand side of \eqref{22} is non-negative. Therefore, setting 
$$d= \sum_{i=4}^r d(\Pi_i)$$ to simplify notation, 
one obtains a solution, for example, by taking $u_1 = k - d - 1$ and assigning all other $u_i$, $v_j$, and $w_e$ the value one, in such a way that \eqref{1} is satisfied. The last requirement can always be fulfilled since \eqref{by} implies that the right-hand sides of the equalities in \eqref{1} are integers and the number $n/a$ satisfies 
$$\frac{n}{a}= \frac{n(a,b,c)}{2a}k=\frac{l(a,b,c)}{a}k\geq k>k-d-1.$$ 
The  branch datum corresponding to this specific solution is   
\begin{equation}
 \l{brd} \Pi(k)=\big((a(k-d-1), a^{n/a-(k-d-1)}), (b^{n/b}), (c^{n/c}), \Pi_4, \dots \Pi_r,  n,0\big).  \qedhere 
\end{equation}
\end{proof}

Concrete examples of non-realizable branch data of the form \eqref{brd} for the triples $(2,3,3)$ and $(2,3,5)$ with $r = k = 3$ are given in the introduction. In addition to these, let us present a similar example
for the triple $(2,3,3)$ with $r=4$, taking
$\Pi_4=(2,1^{n-2})$. In this case, $d = 1$, and \eqref{w2} is satisfied for $k \geq 3$. In particular, for $k = 3$ we obtain the following non-realizable branch datum: 
$$\big((2^9),(3^6),(3^6),(2,1^{16}),18,0\big).$$ 

Note that the set of data for which non-realizability follows from Corollary~\ref{t0} is certainly not restricted to the case where $\deg f$ is not divisible by $\deg \theta_{\f O}$, as in Proposition~\ref{proo}. Let us consider, for instance, the branch datum
\eqref{brd} with $k$ even and
\be \l{w3}
d\leq\frac{k}{2}-2.
\ee
By Corollary~\ref{t0}, a rational function $f$ with this
branch datum would have the form
$f=\theta_{\f O}\circ q$, where $$\deg q=\frac{n}{n(a,b,c)}=\frac{k}{2}.$$
At the point corresponding to the entry $a(k-d-1)$,
the chain rule gives local degree $k-d-1$ for $q$,
since $\theta_{\f O}$ has local degree $a$ over the
corresponding critical value.
However, this is impossible because 
$$k-d-1>k/2=\deg q$$ by \eqref{w3}.
Thus, under the above conditions, \eqref{brd} is
non-realizable. The simplest such example for the triple $(3,3,2)$ is the branch datum
$$\big((9,3^5),(3^8),(2^{12}),24,0\big),$$
where $a=b=3$, $c=2$, $k=4$, and $d=0$.

We conclude this subsection with examples of
non-realizable branch data obtained using
assertion~II.2 of Theorem~\ref{main}. 
Consider the series of branch data
\be \l{thd}
\big((2^{k}),(2^{k-2},1,3),(k,k),2k,0\big),
\qquad k\geq2,
\ee
mentioned in the introduction, whose non-realizability
is known (see \cite{hu}, \cite{pp}).
Assume that a rational function $f$ realizing \eqref{thd}
exists, and let $z_1,z_2,z_3$ be the critical values of $f$
corresponding, respectively, to the three partitions
in \eqref{thd}. Then for the orbifold $\f O$ on $\C\P^1$ defined by the equalities 
$$\nu(z_1)=2, \  \ \ \ \nu(z_2)=2, \  \ \ \ \nu(z_3)=k,$$ the orbifold $f^*(\f O)$ has the signature $\{2,2\}$,
so the equalities in assertion~II.2 of
Theorem~\ref{main} hold. Moreover, it follows from  
$$\deg \theta_{f^*(\f O)}=2, \ \ \ \    \deg \theta_{{\f O}}=2k, \ \ \ \ \deg f=2k$$
that $\deg w=k$ and $\deg t=2.$  

Let $x$ be the point of local degree $3$ in the fiber
$f^{-1}(z_2)$. Since $\deg t=2$, the chain rule implies
that $$\deg_xt=1, \qquad   \deg_{t(x)}w=3.$$
The fiber $t^{-1}(t(x))$ therefore contains another
point $x'\neq x$ at which $t$ is unramified.
Hence $f(x')=z_2$ and $\deg_{x'}f=3$, contradicting the
fact that the partition $(2^{k-2},1,3)$ contains only
one entry equal to $3$.

\subsection{\l{s31} Non-realizability from bad pullbacks}
We now turn to assertion~I of Theorem~\ref{main}.
Choosing a good target orbifold whose pullback would
be bad gives the following non-realizability criteria.

\bc\l{c2} 
A branch datum  $\Pi=(\Pi_1,\Pi_2,\dots, \Pi_r,n,0)$ 
is non-realizable whenever there exist integers  $a,b,c$, each greater than one, such that one of the following conditions holds:

\textnormal{\textsc{Condition 1:}}
\begin{enumerate}
    \item All entries of \(\Pi_1\)  are divisible by \(a\).
    \item All entries of \(\Pi_2\) are divisible by \(b\).
    \item All entries of \(\Pi_3\) except one are  divisible by \(c\).
\end{enumerate}

\textnormal{\textsc{Condition 2:}}
\begin{enumerate}
    \item All entries of \(\Pi_1\) are divisible by \(a\).
    \item All entries of \(\Pi_2\) are divisible by \(b\).
    \item All entries of \(\Pi_3\)  are divisible by \(c\) except for two entries \(u\) and \(v\).
    \item \(\gcd(u, c) \neq \gcd(v, c)\).
\end{enumerate}

\textnormal{\textsc{Condition 3:}}
\begin{enumerate}
    \item All entries of \(\Pi_1\) are divisible by \(a\).
    \item All entries of \(\Pi_2\) are divisible by \(b\) except for one entry  \(u\).
    \item All entries of \(\Pi_3\) are divisible by \(c\)  except for one entry  \(v\).
    \item \(\frac{b}{\gcd(u, b)} \neq \frac{c}{\gcd(v, c)}\).
\end{enumerate}

\ec
\pr
Assume that there exists a rational function $f$ whose
branch datum satisfies one of the above conditions.
Let $\f O$ be the orbifold defined by
$$
\nu(z_1)=a,\qquad \nu(z_2)=b,\qquad \nu(z_3)=c,
$$
where $z_1,z_2,z_3$ are the critical values of $f$
corresponding, respectively, to the partitions
$\Pi_1,\Pi_2,\Pi_3$. 
Then $\f O$ is good, whereas $f^*(\f O)$ is bad,
contrary to assertion~I of Theorem~\ref{main}.
\qed

\vskip 0.2cm

Corollary \ref{c2} allows us to construct many non-realizable branch data. For instance, by slightly modifying the series \eqref{iz0}, we obtain the non-realizable series 
$$\big((2^k),(2^{k-2},4),(k-s, k-s, 2s),2k,0\big), \ \ \  k> 3s,$$ 
by applying Condition~1 with $a=b=2$ and $c=k-s.$

Similarly, we see that, along with the series \eqref{iz}, the series  \be \l{z1} \big((2^l), (3^k,s), (5^m), n,0\big), \ \ \ s\not\equiv 0 \ (\mod 3), \ee
 \be \l{z2} \big((2^l,s), (3^k), (5^m), n,0\big), \ \ \ s\not\equiv 0 \ (\mod 2), \ee
or, for instance, the series 
\be \l{z3} \big((2^l), (3^k,d), (5^m,s), n,0\big),\ee where at least one of the conditions $d\not\equiv 0 \ (\mod 3)$, $s\not\equiv 0 \ (\mod 5)$ holds, are also  non-realizable. 

Many other examples can be obtained through calculations similar to those performed in Section \ref{s32}, by 
 considering branch data such that in the first three partitions all entries but one or two are divisible by suitable numbers $a,b,c$, though not necessarily equal to those numbers as in \eqref{iz}, \eqref{z1}, \eqref{z2}, \eqref{z3}.   

\subsection{\l{s33} A decomposability criterion} 
The factorization in Theorem~\ref{main}  yields
the following decomposability criterion mentioned in the introduction. 

\bc \l{co} Let $f$ be a rational function. Suppose that for some  good orbifold $\f O$ the conditions $\chi(f^*(\f O))>0$ and  $\nu(f^*(\f O))\neq \nu(\f O)$ hold. Then $f$ is either decomposable or a compositional
left factor of $\theta_{\f O}$.
\ec 
\pr
If $f$ is indecomposable, the first equality in
assertion~II.2 of Theorem~\ref{main} implies that either
$\deg w=1$ or $\deg t=1$.
In the first case, 
equality \eqref{noth} holds by Lemma~\ref{lmo},
contradicting the assumption.
Therefore, $\deg t=1$, which implies that $f$ is
a compositional left factor of $\theta_{\f O}$
by the second equality in the same assertion.
\qed

Corollary~\ref{co} unifies and extends several decomposability criteria previously obtained in the study of monodromy groups of indecomposable rational functions (see \cite{gur}, Proposition~2.0.16--Corollary~2.0.18, and \cite{nz}, Lemma~9.1). The cited criteria concern particular ramification patterns, expressed through divisibility conditions on local degrees and restrictions on exceptional cycles rather than in terms of orbifolds. The language of orbifolds provides a particularly convenient framework for treating these questions uniformly. Specifically,  Corollary~\ref{co} encompasses all spherical target signatures and, for each target, all pullback signatures of positive Euler characteristic distinct from the target signature. Its scope therefore extends both to additional target signatures and to further ramification patterns over previously considered targets.

For example, writing $\f O_1=f^*(\f O_2)$, the three
cases of Lemma~9.1 in \cite{nz}
(which is an explicit version of the translation method used in
\cite{gur}, Proposition~2.0.16 and Corollary~2.0.17)
follow from Corollary~\ref{co} for the following choices:
\begin{itemize}
\item $\nu(\f O_2)=\{p,p\}$, where $p$ is a prime,
and $\f O_1$ is non-ramified.
\item $\nu(\f O_2)=\{p,2,2\}$, where $p$ is a prime,
and $\nu(\f O_1)=\{2,2\}$.
\item $\nu(\f O_2)=\{2,3,3\}$ and
$\nu(\f O_1)=\{3,3\}$.
\end{itemize}
In particular, the octahedral and icosahedral target
signatures, $\{2,3,4\}$ and $\{2,3,5\}$, do not occur
among these three cases, whereas both are covered
by Corollary~\ref{co}.

Corollary~2.0.18 of \cite{gur} applies to the
icosahedral target $\nu(\f O_2)=\{2,3,5\}$, but imposes
two additional restrictions. First, the three
corresponding branch permutations must have orders
$2,3,5$, respectively. Second, the sum of the numbers
of fixed points of these permutations must be at most
two, which, under the first restriction, is equivalent
to $\f O_1$ having at most two singular points. 
Moreover, as formulated, Corollary~2.0.18 forces decomposability
only for functions with more than three critical  values.

For example, if a rational function $f$ realizes the
branch datum
$$
\big((4,2^{15},1),(3^{11},1^2),(5^7),35,0\big),
$$
then $f$ has exactly three critical values, and both
restrictions above fail since the corresponding branch
permutations have orders $4,3,5$ and a total of three fixed points. Nevertheless, assigning weights $2,3,5$, respectively,
to the critical values corresponding to the three
partitions above gives an orbifold $\f O_2$ for which
$\nu(f^*(\f O_2))=\{2,3,3\}$. 
Corollary~\ref{co} therefore implies that $f$ is
decomposable, since its degree $35$ does not divide
$\deg\theta_{\f O_2}=60$, ruling out the exceptional case.

This branch datum is indeed realized by a composition
$f=w\circ t$, where $w$ is a degree-five orbifold covering
with source and target signatures $\{2,3,3\}$ and
$\{2,3,5\}$, respectively, with branch datum
$$
\big((2,2,1),(3,1,1),(5),5,0\big),
$$
and $t$ is a rational function with branch datum
$$
\big((4,2,1),(3,3,1),(3,3,1),7,0\big).
$$
The three partitions in the branch datum of $t$
correspond to the singular points of the intermediate
orbifold of weights $2,3,3$, respectively.

\section{\l{s4} Proof and applications of Theorem~\ref{t3}}
\subsection{\l{s41} Proof of Theorem~\ref{t3}} 
Let $\{S,p,q\}$ be a component of the fiber product of $f$ and
$g=\theta_{\f O}$. 
Since $f^*(\f O)$ is non-ramified, Theorem~\ref{t5} shows that $p$
is unramified. Then the Riemann--Hurwitz formula gives $g(S)=1$. Since an unramified
covering between compact Riemann surfaces of genus one
is Galois with an Abelian deck transformation group  (see \cite{silv}, Theorem~III.4.10), we conclude that $p$ is a Galois covering and
$H=\Aut(S,p)$ is Abelian.

By Theorem~\ref{gc}, there is an injective homomorphism
$$
\phi:H\to G=\Aut(\C\P^1,\theta_{\f O})
$$
satisfying
\be \l{sa}
q\circ\sigma=\phi(\sigma)\circ q
\ee
for every $\sigma\in H$.
Set $K=\phi(H)$ and let
$v:\C\P^1\to\C\P^1/K$ be the quotient map.
As in the proof of Theorem \ref{main}, identifying $\C\P^1/K$ with $\C\P^1$, we have
$\theta_{\f O}=w\circ v$ for a rational function $w$.
Moreover, since \eqref{sa} implies that $v\circ q$ is $H$-invariant, the equality 
$$v\circ q=t\circ p$$ holds for a holomorphic map $t:R\to\C\P^1$.
Consequently, $$f\circ p=\theta_{\f O}\circ q=w\circ t\circ p,$$ and hence
$f=w\circ t$.

Set $\f O'=w^*(\f O)$. By Corollary~\ref{indu2}, both
$v:\C\P^1\to\f O'$ and $w:\f O'\to\f O$ are covering maps
between orbifolds. Thus $v$ is a universal covering of $\f O'$.
Choosing $\theta_{w^*(\f O)}=v$, we obtain
\be\l{bur}
f=w\circ t,\qquad
\theta_{\f O}=w\circ\theta_{w^*(\f O)}.
\ee
Finally, since $K$ is a finite Abelian group of
M\"obius transformations, it is cyclic or isomorphic
to the Klein four-group. Therefore,
$\nu(w^*(\f O))$ is $\{d,d\}$ for some $d\geq1$,
or $\{2,2,2\}$. \qed

\subsection{\l{s42} Applications to non-realizability and decomposability}
As an example illustrating Theorem \ref{t3}, we show that the series of branch data  
\be \l{wbd} \big((2^{k},k+3),(3^{k+1}),(3^{k+1}),3k+3,1\big), \ \ \ k\equiv 1\, ({\rm mod}\, 4),\ee is non-realizable. 

Assume that $f$ is a holomorphic map realizing \eqref{wbd},
and let $z_1,z_2,z_3$ be the critical values of $f$
corresponding, respectively, to the three partitions
in \eqref{wbd}. Then for the orbifold $\f O$ on $\C\P^1$ defined by the equalities 
$$\nu(z_1)=2, \  \ \ \ \nu(z_2)=3, \  \ \ \ \nu(z_3)=3,$$ the orbifold $f^*(\f O)$ is non-ramified, implying by Theorem \ref{t3} that equalities \eqref{bur} hold. 
On the other hand, 
any decomposition of $\theta_{\f O}$ into indecomposable
rational functions of degree at least two has either
the form
$$
\theta_{\f O}=w_1\circ w_2\circ w_3,
$$
where $\deg w_1=3$, $\deg w_2=2$, and $\deg w_3=2$,
or the form
$
\theta_{\f O}=s_1\circ s_2,
$
where $\deg s_1=4$ and $\deg s_2=3$.
Moreover, the branch datum of $w_1$ is
$\big((3),(3),3,0\big)$  (see e.g. \cite{gen},  Section 4.3).

Since $k\equiv1\pmod4$, the degree of $f$ is not divisible
by four. Therefore, 
$f$ admits a decomposition $f=w\circ t$, where $w$ has
degree three and branch datum $\big((3),(3),3,0\big)$,
while $\deg t=k+1$. By the chain rule, $z_1$ cannot be
a critical value of $w$, since the partition over $z_1$
contains an entry equal to two. Hence $t$ must have
a point of local degree $k+3$, which is impossible
because $$k+3>k+1=\deg t.$$

As another example of using  Theorem \ref{t3} for proving non-realizability, we
consider the series of branch data  
\be \l{koro} \big((2^{3k+6}),(3^{2k+4}),(3,9,6^{k}),6k+12,1\big), \ \ \ k\equiv 1\, ({\rm mod}\, 2).\ee Arguing as above one can see that if $f$ is a holomorphic map  realizing \eqref{koro}, then $f=w\circ t$, where $w$ is a rational function of degree 3 with the branch datum $\big((3),(3),3,0\big)$ and $t$ is a holomorphic map of degree $2k+4.$ Moreover, if  
 $z_1$ is a critical value of $f$ corresponding to the partition $(2^{3k+6})$, then $z_1$ is not a critical value of $w$. 
 
 Over each of the three points of $w^{-1}(z_1)$,
the map $t$ has partition $(2^{k+2})$, by the chain
rule and the equality $\deg t=2k+4$.
Hence, the two critical values of $w$ must be the critical
values of $f$ corresponding to the second and third
partitions in \eqref{koro}. Since $w$ is totally
ramified over these values, the corresponding
partitions for $t$ are obtained by dividing all
entries by three, giving $(1^{2k+4})$ and
$(1,3,2^k)$, respectively.
The first of these partitions is trivial, and $t$
has no other critical values. Thus, its branch datum is
$$
\big((2^{k+2}),(2^{k+2}),(2^{k+2}),
(1,3,2^k),2k+4,1\big).
$$
However, this branch datum is known to be non-realizable
(see \cite{ko}, Section~5). Therefore, the branch data
\eqref{koro} are non-realizable as well.

Extending the notion of decomposability to holomorphic
maps between compact Riemann surfaces, we obtain the
following consequence of Theorem~\ref{t3}, analogous
to Corollary~\ref{co}.

\bc  Let $R$ be a compact Riemann surface of genus one and \linebreak  $f:R\rightarrow \C\P^1$ a holomorphic map. Assume that for some orbifold  $\f O$   on $\C\P^1$ with the signature 
 $\{2,2,d\}$, $d>2,$ $\{2,3,3\},$ $\{2,3,4\},$ or $\{2,3,5\}$ the orbifold $f^*(\f O)$ is non-ramified.  Then $f$ is decomposable. 
\ec
\pr If $f$ is indecomposable, then in the first equality of \eqref{bur}, either $\deg w = 1$ or $\deg t = 1$. Both of these assumptions lead to a contradiction. Indeed, if $\deg w = 1$, then the second equality in \eqref{bur} implies, by Lemma \ref{lmo}, that the signature of $\f O$ is $\{d,d\}$, with $d \geq 1$, or $\{2,2,2\}$, which contradicts the assumption. The equality $\deg t = 1$ is also impossible, since $w$ is defined on the Riemann sphere, while $f$ is defined on a torus. \qed

\section{\l{s5} Halphen's theorem}
In this section, we deduce from Theorem~\ref{main}
Halphen's theorem (see \cite{hal} or \cite{arn})
concerning polynomial solutions of the generalized
Fermat equation \be \l{eq} X^a + Y^b = Z^c,\ee 
where  $(a,b,c)$ is a triple of integers  greater than one.

We begin by constructing solutions of \eqref{eq} from rational Galois coverings. Assume that $(a,b,c)$ satisfies  $\chi(a,b,c)>0$, and let $\theta_{\f O}$ be a universal covering of an orbifold $\f O$ defined by the equalities
\be \l{kro} \nu(1)=a,  \ \ \ \ \nu(\infty)=b, \ \ \ \ \nu(0)=c.\ee 
Since $\theta_{\f O}$ is uniform, there exist coprime polynomials 
$P$, $Q$, $R$  
such that 
\be \l{kro1} \theta_{\f O}=\frac{R^c}{P^b}, \ \ \  \theta_{\f O}-1=\frac{Q^a}{P^b}.\ee

It is clear that 
$$Q^a(z)+P^b(z)=R^c(z).$$ Moreover, this equality remains true after the substitution $z=U/V$, where $U$ and $V$ are coprime polynomials. This implies that for any non-zero complex numbers $\alpha,\beta,\gamma$ such that \be \l{abc} \alpha^a=\beta^b=\gamma^c,\ee the rational functions 
\be \l{kro3} 
X = \alpha V^{n/a}Q\left(\frac{U}{V}\right),
\quad
Y = \beta V^{n/b}P\left(\frac{U}{V}\right),
\quad
Z = \gamma V^{n/c}R\left(\frac{U}{V}\right),
\ee
where $n=\deg\theta_{\f O}$, are polynomials satisfying
\eqref{eq}. Moreover, these polynomials are pairwise
coprime. Indeed, the weighted degrees $a\,\deg Q$,
$b\,\deg P$, and $c\,\deg R$ do not exceed $n$,
and at least two of them equal $n$. Hence, since $U$ and $V$
are coprime, at a zero of $V$ at most one of these
polynomials vanishes. Away from these zeros, their
pairwise coprimality follows from that of $P,Q,R$.

In the above notation, Halphen's theorem takes
the following form.

\bt Let $(a, b, c)$ be a triple of integers greater than one. Then equation \eqref{eq} has no solutions in coprime non-constant polynomials $X, Y, Z$ unless $\chi(a, b, c) > 0$. On the other hand, if $\chi(a, b, c) > 0$, then every such solution of \eqref{eq} has the form \eqref{kro3}, 
where  $(P,Q,R)$ is any  triple of coprime polynomials satisfying \eqref{kro1} 
 for the orbifold $\f O$ defined by \eqref{kro}.
\et
\pr
Let $X,Y,Z$ be coprime non-constant polynomials satisfying
\eqref{eq}, and let $\f O$ be the orbifold defined
by \eqref{kro}. Since $\f O$ has three singular points,
it is good. Set
\be \l{say1}
F=\frac{Z^c}{Y^b}.
\ee
Then
\be \l{say2}
F-1=\frac{X^a}{Y^b}.
\ee
Since $X,Y,Z$ are pairwise coprime, the local degrees
of $F$ at finite points over $1$, $\infty$, and $0$
are divisible by $a$, $b$, and $c$, respectively.
Thus $F^*(\f O)$ is either non-ramified or has
$\infty$ as its only singular point.
The latter case is impossible by assertion~I
of Theorem~\ref{main}. Therefore, 
$\chi(\f O)>0$ and $F=\theta_{\f O}\circ q$ for some rational function $q$ by Corollary \ref{t0}. 

Writing the non-constant rational function $q$ in reduced form
as $q=U/V$, we obtain  
$$\frac{Z^c}{Y^b}=F=\theta_{\f O} \left(\frac{U}{V}\right)= \frac{R^c\left(\frac{U}{V}\right)}{P^b\left(\frac{U}{V}\right)}=\frac{R^c\left(\frac{U}{V}\right)V^{n}}{P^b\left(\frac{U}{V}\right)V^{n}}=\frac{\left(R\left(\frac{U}{V}\right)V^{n/c}\right)^c}{\left(P\left(\frac{U}{V}\right)V^{n/b}\right)^b}.$$ 
Since  $P,Q,R$ and $X,Y,Z$ are triples of coprime polynomials, this   implies that 
\be \l{mo} Y = \beta V^{n/b}P\left(\frac{U}{V}\right), \ \ \ 
Z = \gamma V^{n/c}R\left(\frac{U}{V}\right),\ee where $\beta^b=\gamma^c.$
Similarly, 
$$\frac{X^a}{Y^b}=F-1= \frac{Q^a\left(\frac{U}{V}\right)}{P^b\left(\frac{U}{V}\right)}=\frac{Q^a\left(\frac{U}{V}\right)V^{n}}{P^b\left(\frac{U}{V}\right)V^{n}}=\frac{\left(Q\left(\frac{U}{V}\right)V^{n/a}\right)^a}{\left(P\left(\frac{U}{V}\right)V^{n/b}\right)^b},$$ implying that 
\be \l{mo1} Y = \t\beta V^{n/b}P\left(\frac{U}{V}\right), \ \ \ 
X = \alpha V^{n/a}Q\left(\frac{U}{V}\right),\ee
where $\t\beta^b=\alpha^a.$ 
Finally, the first equalities in \eqref{mo} and \eqref{mo1} imply  that $\t\beta=\beta.$ \qed
\vskip 0.2cm
Note that the non-existence assertion in Halphen's
theorem follows easily from the Mason--Stothers theorem
(see \cite{mas}). In contrast, the parametrization
\eqref{kro3}, which is also part of Halphen's classical
theorem, requires additional arguments (see \cite{arn},
Section~28.2.3, Theorem~28.2.1 and the discussion
following its proof).

Finally, note that   the Halphen theorem   implies in turn  Corollary \ref{t0}. 
Indeed, let $f$ satisfy the assumptions of Corollary~\ref{t0}.
By Proposition~\ref{p1}, we have $\chi(\f O)>0$.
The cases where $\f O$ is non-ramified or has two
singular points are immediate, so it remains to
consider the case of three singular points.  Assuming without loss of generality that the orbifold $\f O$ is defined by the equalities \eqref{kro} and that
$\infty$ is not a critical point of $f$, we see that  there exist coprime non-constant  polynomials 
$X,Y,Z$ such that equalities \eqref{say1} and \eqref{say2} hold with $F=f$. Thus, $(X,Y,Z)$ is a solution of \eqref{eq}, implying
by the Halphen theorem that
$$f=Z^c/Y^b=
\frac{\left(\gamma V^{n/c}R\left(\frac{U}{V}\right)\right)^c}
{\left(\beta V^{n/b}P\left(\frac{U}{V}\right)\right)^b}
=\frac{R^c\left(\frac{U}{V}\right)}{P^b\left(\frac{U}{V}\right)}
=\theta_{\f O}\left(\frac{U}{V}\right).
$$

\section*{Acknowledgements}
The author is grateful to Nikolai Adrianov for drawing
attention to \cite{adr}, and to Danny Neftin for pointing
out references \cite{gur} and \cite{nz}. He also thanks Hao Zheng for providing
computational results related to his paper \cite{zh} that were
not included in the published version.

\end{document}